# Shape restricted regression with random Bernstein polynomials

### I-Shou Chang[1,*], Li-Chu Chien[2,*], Chao A. Hsiung[2], Chi-Chung Wen[3] and Yuh-Jenn Wu[4]

*National Health Research Institutes, Tamkang University, and Chung Yuan Christian University, Taiwan*

**Abstract:** Shape restricted regressions, including isotonic regression and concave regression as special cases, are studied using priors on Bernstein polynomials and Markov chain Monte Carlo methods. These priors have large supports, select only smooth functions, can easily incorporate geometric information into the prior, and can be generated without computational difficulty. Algorithms generating priors and posteriors are proposed, and simulation studies are conducted to illustrate the performance of this approach. Comparisons with the density-regression method of Dette et al. (2006) are included.

## 1. Introduction

Estimation of a regression function with shape restriction is of considerable interest in many practical applications. Typical examples include the study of dose response experiments in medicine and the study of utility functions, product functions, profit functions and cost functions in economics, among others. Starting from the classic works of Brunk [4] and Hildreth [17], there exists a large literature on the problem of estimating monotone, concave or convex regression functions. Because some of these estimates are not smooth, much effort has been devoted to the search of a simple, smooth and efficient estimate of a shape restricted regression function. Major approaches to this problem include the projection methods for constrained smoothing, which are discussed in Mammen et al. [24] and contain smoothing splines methods and kernel and local polynomial methods and others as special cases, the isotonic regression approach studied by Mukerjee [25], Mammen [22, 23] and others, the tilting method proposed by Hall and Huang [16], and the density-regression method proposed by Dette, Neumeyer and Pilz [11]. We note that both of the last two methods enjoy the same level of smoothness as the unconstrained counterpart and are applicable to general smoothing methods. Besides, the density-regression method is particularly computationally efficient and has a wide applicability. In

*Partially supported by NSC grant, NSC 94-3112-B-400-002-Y.

[1]Institute of Cancer Research and Division of Biostatistics and Bioinformatics, National Health Research Institutes, 35 Keyan Road, Zhunan Town, Miaoli County 350, Taiwan, e-mail: ischang@nhri.org.tw

[2]Division of Biostatistics and Bioinformatics, National Health Research Institutes, 35 Keyan Road, Zhunan Town, Miaoli County 350, Taiwan, e-mail: lcchien@nhri.org.tw; hsiung@nhri.org.tw

[3]Department of Mathematics, Tamkang University, 151 Ying-chuan Road, Tamsui Town, Taipei County, 251, Taiwan, e-mail: ccwen@mail.tku.edu.tw

[4]Department of Applied Mathematics, Chung Yuan Christian University, 200 Chung Pei Road, Chung Li City 320, Taiwan, e-mail: yuhjenn@cycu.edu.tw

fact, the density-regression method was used to provide an efficient and smooth convex estimate for convex regression by Birke and Dette [3].

This paper studies a nonparametric Bayes approach to shape restricted regression where the prior is introduced by Bernstein polynomials. This prior features the properties that it has a large support, selects only smooth functions, can easily incorporate geometric information into the prior, and can be generated without computational difficulty. We note that Lavine and Mockus [18] also discussed Bayes methods for isotonic regression, but the prior they use is a Dirichlet process, whose sample paths are step functions. In addition to the above desirable properties, our approach can be applied to quite general shape restricted statistical inference problems, although we consider mainly isotonic regression and concave (convex) regression in this paper. To facilitate the discussion, we first introduce some notations as follows.

For integers $0 \le i \le n$, let $\varphi_{i,n}(t) = C_i^n t^i (1-t)^{n-i}$, where $C_i^n = n!/(i!(n-i)!)$. $\{\varphi_{i,n} \mid i = 0, \ldots, n\}$ is called the Bernstein basis for polynomials of order $n$. Let $B_n = \mathbb{R}^{n+1}$ and $\mathcal{B} = \bigcup_{n=1}^{\infty} (\{n\} \times B_n)$. Let $\pi$ be a probability measure on $\mathcal{B}$. For $\tau > 0$, we define $\mathbf{F} : \mathcal{B} \times [0, \tau] \longrightarrow \mathbb{R}^1$ by

$$\mathbf{F}(n, b_{0,n}, \ldots, b_{n,n}, t) = \sum_{i=0}^{n} b_{i,n} \varphi_{i,n}(\frac{t}{\tau}), \qquad (1.1)$$

where $(n, b_{0,n}, \ldots, b_{n,n}) \in \mathcal{B}$ and $t \in [0, \tau]$. We also denote (1.1) by $F_{b_n}(t)$ if $b_n = (b_{0,n}, \ldots, b_{n,n})$. The probability measure $\pi$ is called a Bernstein prior, and $\mathbf{F}$ is called the random Bernstein polynomial for $\pi$. It is a stochastic process on $[0, \tau]$ with smooth sample paths. Important references for Bernstein polynomials include Lorentz [20] and Altomare and Campiti [2], among others. It is well-known that Bernstein polynomials are popular in curve design (Prautzsch et al. [28]).

Bernstein basis have played important roles in nonparametric curve estimation and in Bayesian statistical theory. Good examples of the former include Tenbusch [31] and some of the references therein. For the latter, we note that Beta density $\beta(x; a, b) = x^{a-1}(1-x)^{b-1}/B(a,b)$ is itself a Bernstein polynomial and mixtures of Beta densities of the kind $\sum_{j=1}^{n} w_j \beta(x; a_j, b_j)$, where $w_j \ge 0$ are random with $\sum_{j=1}^{n} w_j = 1$, were used to introduce priors that only select smooth density functions on $[0, 1]$; see Dalal and Hall [8], Diaconis and Ylvisaker [9], and Mallik and Gelfand [21], and references therein. We also note that Petrone [26] and Petrone and Wasserman [27] studied priors on the set of distribution functions on $[0, 1]$ that are specified by $\sum_{i=0}^{n} \tilde{F}(i/n) \varphi_{i,n}(t)$ with $\tilde{F}$ being a Dirichlet process; this prior was referred to as a Bernstein-Dirichlet prior.

For a continuous function $F$ on $[0, \tau]$, $\sum_{i=0}^{n} F(i\tau/n) \varphi_{i,n}(t/\tau)$ is called the $n$-th order Bernstein polynomial of $F$ on $[0, \tau]$. We will see in Section 2 that much of the geometry of $F$ is preserved by its Bernstein polynomials and very much of the geometry of a Bernstein polynomial can be read off from its coefficients. This together with Bernstein-Weierstrass approximation theorem suggests the possibility of a Bernstein prior on a space of functions with large enough support and specific geometric properties. These ideas were developed in Chang et al. [7] for Bayesian inference of a convex cumulative hazard with right censored data.

The purpose of this paper is to indicate that the above ideas are also useful in the estimation of shape restricted regressions. The regression model we consider in this paper assumes that given $F_{b_n}$ satisfying certain shape restriction,

$$Y_{jk} = F_{b_n}(X_k) + \epsilon_{jk},$$

where $X_k$ are design points, $Y_{jk}$ are response variables and $\epsilon_{jk}$ are errors.

We will investigate isotonic regression and concave (convex) regression in some detail. In particular, we will examine the numerical performance of this approach and offer comparison with the newly-developed density-regression method of Dette et al. [11] in simulation studies. We also indicate without elaboration that this approach can also be used to study regression function that is smooth and unimodal or that is smooth and is constantly zero initially, then increasing for a while, and finally decreasing. We note that the latter can be used to model the time course expression level of a virus gene, as discussed in Chang et al. [5]; a virus gene typically starts to express after it gets into a cell and its expression level is zero initially, increases for a while and then decreases.

That it is very easy both to generate the prior and the posterior for inference with Bernstein prior is an important merit of this approach. Because the prior is defined on the union of subspaces of different dimension, we adapt Metropolis-Hastings reversible jump algorithm (MHRA) (Green [15]) to calculate the posterior distributions in this paper.

This paper is organized as follows. Section 2 introduces the model, provides statements that exemplify the relationship between the shape of the graph of (1.1) and its coefficients $b_0, \ldots, b_n$, and gives conditions under which the prior has full support with desired geometric peoperties. Section 3 illustrates the use of Bernstein priors in conducting Bayesian inference; in particular, suitable Bernstein priors are introduced for isotonic regression and unimodal concave (convex) regression, and Markov chain Monte Carlo approaches to generate posterior are proposed. Section 4 provides simulation studies to compare our methods with the density-regression method. Section 5 is a discussion on possible extensions.

## 2. The Bernstein priors

### 2.1. Bernstein polynomial geometry

Let $F_a(t) = \sum_{i=0}^{n} a_i \varphi_{i,n}(t/\tau)$. This subsection presents a list of statements concerning the relationship between the shape of $F_a$ and its coefficients $a_0, \ldots, a_n$. This list extends that in Chang et al. [7] and is by no means complete; similar statements can be made for a monotone and convex or a monotone and sigmoidal function, for example. They are useful in taking into account the geometric prior information for regression analysis.

**Proposition 1.** *(i) (Monotone) If $a_0 \leq a_1 \leq \cdots \leq a_n$, then $F_a'(t) \geq 0$ for every $t \in [0, \tau]$.*

*(ii) (Unimodal Concave) Let $n \geq 2$. If $a_1 - a_0 > 0, a_n - a_{n-1} < 0$, and $a_{i+1} + a_{i-1} \leq 2a_i$, for every $i = 1, \ldots, n-1$, then $F_a'(0) > 0, F_a'(\tau) < 0$, and $F_a''(t) \leq 0$ for every $t \in [0, \tau]$.*

*(iii) (Unimodal) Let $n \geq 3$. If $a_0 = a_1 = \cdots = a_{l_1} < a_{l_1+1} \leq a_{l_1+2} \leq \cdots \leq a_{l_2}$ and $a_{l_2} \geq a_{l_2+1} \geq \cdots \geq a_{l_3} > a_{l_3+1} = \cdots = a_n$ for some $0 \leq l_1 < l_2 < l_3 \leq n$, then there exists $s \in (0, \tau)$ such that $s$ is the unique maximum point of $F_a$ and $F_a$ is strictly increasing on $[0, s]$ and strictly decreasing on $[s, \tau]$. Furthermore, if $l_1 > 0$, then $F_a^{(i)}(0) = 0$ for $i = 1, \ldots, l_1$, and if $l_3 < n - 1$, then $F_a^{(i)}(\tau) = 0$ for $i = 1, \ldots, n - l_3 - 1$.*

In this paper, derivatives at 0 and $\tau$ are meant to be one-sided. We note that, in Proposition 1, (ii) provides a sufficient condition under which $F_a$ is a unimodal concave function and (iii) provides a sufficient condition under which $F_a$ is a unimodal smooth function whose function and derivative values at 0 and $\tau$ may be prescribed to be 0. Let $F_a$ satisfy (iii) with $\tau = 1$ and $l_1 > 1$, and let $\overline{F}$ be defined by $\overline{F}(t) = F_a(t - \tau_1/\tau - \tau_1)I_{(\tau_1, \tau)}(t)$ for some $\tau_1 \in (0, \tau)$, then $\overline{F}$ is a non-negative smooth function on $[0, \tau]$ and it is zero initially for a while, then increases, and finally decreases. As we mentioned earlier, functions with this kind of shape restriction are useful in the study of time course expression profile of a virus gene. We note that the results presented here are for concave regressions or unimodal regressions and similar results hold for convex regressions.

*Proof.* Without loss of generality, we assume $\tau = 1$. Noting that the derivatives $F_a'(t) = n\sum_{i=0}^{n-1}(a_{i+1} - a_i)\varphi_{i,n-1}(t)$ and $F_a''(t) = n(n-1)\sum_{i=0}^{n-2}(a_{i+2} - 2a_{i+1} + a_i)\varphi_{i,n-2}(t)$, we obtain (i) and (ii) immediately. We now prove (iii) for the case $l_1 = 1$ and $l_3 = n - 1$; the proofs for other cases are similar and hence omitted.

Let $\varphi(x) = \sum_{i=0}^{n-1} C_i^{n-1}(a_{i+1} - a_i)x^i$, then $F_a'(t) = n(1-t)^{n-1}\varphi(\frac{t}{1-t})$. Using $a_0 = a_1 < a_2 \leq \cdots \leq a_l$ and $a_l \geq a_{l+1} \geq \cdots \geq a_{n-1} > a_n$, we know $\varphi$ is a polynomial whose number of sign changes is exactly one. This together with Descartes' sign rule (Anderson et al. [1]) implies that $\varphi$ has at most one root in $(0, \infty)$, and hence, $F_a'$ has at most one root in $(0, 1)$.

Using $a_0 = a_1$ and $a_2 - a_1 > 0$, we know $F_a'(\epsilon) > 0$ for $\epsilon$ being positive and small enough. This combined with $F_a'(1) = n(a_n - a_{n-1}) < 0$ shows that $F_a'$ has at least one root in $(0, 1)$. Thus, $F_a'$ has exactly one root $s$ in $(0, 1)$, and $F_a'$ is positive on $(0, s)$ and negative on $(s, 1]$. Therefore, the conclusion of (iii) follows. This completes the proof. $\square$

The following proposition complements Proposition 1 and provides Bernstein-Weierstrass approximations for functions with specific shape restrictions.

Let $I_n^{(1)} = \{F_a \mid a \in B_n, a_0 \leq a_1 \leq \cdots \leq a_n\}$, $I_n^{(2)} = \{F_a \mid a \in B_n, a_1 - a_0 > 0, a_n - a_{n-1} < 0, a_{i+1} + a_{i-1} \leq 2a_i, \text{for } i = 1, \ldots, n-1\}$, and $I_n^{(3)} = \{F_a \mid a \in B_n, a_0 = a_1 < a_2 \leq \cdots \leq a_l, a_l \geq a_{l+1} \geq \cdots \geq a_{n-1} > a_n, \text{for some } l = 2, \ldots, n-1\}$. Then we have

**Proposition 2.** *(i) (Monotone) Let $\mathcal{D}_1$ consist of linear combinations of elements in $\bigcup_{n=1}^{\infty} I_n^{(1)}$, with non-negative coefficients. Then the closure of $\mathcal{D}_1$ in uniform norm is precisely the set of increasing and continuous functions on $[0, \tau]$.*

*(ii) (Unimodal Concave) Let $\mathcal{D}_2$ consist of linear combinations of elements in $\bigcup_{n=2}^{\infty} I_n^{(2)}$, with non-negative coefficients. Let $S$ denote the set of all continuously differentiable real-valued functions $F$ defined on $[0, \tau]$ with $F'(0) \geq 0$, $F'(\tau) \leq 0$, and $F'$ decreasing. For two continuously differentiable functions $f$ and $g$, define $d(f, g) = \|f - g\|_\infty + \|f' - g'\|_\infty$, where $\|\cdot\|_\infty$ is the sup-norm for functions on $[0, \tau]$. Then the closure of $\mathcal{D}_2$, under $d$, is $S$.*

*(iii) (Unimodal) Let $\mathcal{D}_3 = \bigcup_{n=3}^{\infty} I_n^{(3)}$. Let $S$ denote the set of all continuously differentiable real-valued functions $F$ defined on $[0, \tau]$ satisfying the properties that $F'(0) = 0$ and that there exists $s \in [0, \tau]$ such that $F'(s) = 0$, $F'(x) \geq 0$ for $x \in [0, s]$, and $F'(x) \leq 0$ for $x \in [s, \tau]$. For two continuously differentiable functions $f$ and $g$, define $d(f, g) = \|f - g\|_\infty + \|f' - g'\|_\infty$, where $\|\cdot\|_\infty$ is the sup-norm for functions on $[0, \tau]$. Then the closure of $\mathcal{D}_3$, under $d$, is $S$.*

*Proof.* We give the proofs for (i) and (ii), and omit the proof for (iii), because it is similar to that for (ii).

(i) It is obvious to see that the closure of $\mathcal{D}_1$ is contained in the set of increasing and continuous functions. We now prove the converse. Let $F$ be increasing and continuous. Taking $a_i = F(i\tau/k)$, for $i = 0, 1, \ldots, k$, we set $F_{(k)}(t) = \sum_{i=0}^{k} a_i \varphi_{i,k}(t/\tau)$. Then $F_{(k)}$ is in $\mathcal{D}_1$. It follows from Bernstein-Weierstrass approximation theorem that $F_{(k)}$ converges to $F$ uniformly. This completes the proof.

(ii) It follows from (ii) in Proposition 1 that $\mathcal{D}_2 \subset S$. Because $S$ is complete relative to $d$, we know $\overline{\mathcal{D}_2} \subset S$. Here $\overline{\mathcal{D}_2}$ is the closure of $\mathcal{D}_2$. We now prove $S \subset \overline{\mathcal{D}_2}$.

Let $F \in S$. Let $d_0 = F'(0) + 1/n^3$, $d_{n-1} = F'(\tau) - 1/n^3$, $d_i = F'(i\tau/(n-1))$, for $i = 1, \ldots, n-2$, and $H_{1,n}(t) = \sum_{i=0}^{n-1} d_i \varphi_{i,n-1}(t/\tau)$. Note that $H_{1,n}(0) > 0$, $H_{1,n}(\tau) < 0$. Using Bernstein-Weierstrass Theorem, we know $H_{1,n}$ converges to $F'$ uniformly.

Let $a_0 = F(0)n/\tau$, $a_i = a_0 + d_0 + \cdots + d_{i-1}$, $H_{0,n}(t) = (\tau/n)\sum_{i=0}^{n} a_i \varphi_{i,n}(t/\tau)$. Then $H'_{0,n}(t) = H_{1,n}(t)$ and $H_{0,n}(0) = F(0)$. Thus $H_{0,n}$ converges uniformly to $F$ (See, for example, Theorem 7.17 in Rudin [30]).

Using the fact that $d_i$ is a decreasing sequence with $d_0 > 0, d_{n-1} < 0$, and

$$a_{i+2} - 2a_{i+1} + a_i = d_{i+1} - d_i,$$

we know $a_1 - a_0 > 0$, $a_n - a_{n-1} < 0$ and $a_{i+2} + a_i \leq 2a_{i+1}$. Thus $H_{0,n}$ is in $\mathcal{D}_2$. This shows that $F$ is in $\overline{\mathcal{D}_2}$. This completes the proof. $\qquad\square$

## 2.2. Bayesian regression

We now describe a Bayesian regression model with the prior distribution, on the regression functions, introduced by random Bernstein polynomials (1.1).

Assume that on a probability space $(\mathcal{B} \times \mathbb{R}^\infty, \mathcal{F}, \mathcal{P})$, there are random variables $\{Y_{jk} \mid j = 1, \ldots, m_k; k = 1, \ldots, K\}$ satisfying the property that, conditional on $\mathbf{B} = (n, b_n)$,

$$Y_{jk} = F_{b_n}(X_k) + \epsilon_{jk}, \tag{2.1}$$

with $\{\epsilon_{jk} \mid j = 1, \ldots, m_k; k = 1, \ldots, K\}$ being independent random variables, $\epsilon_{jk}$ having known density $g_k$ for $j = 1, \ldots, m_k$, $\mathbf{B}$ being the projection from $\mathcal{B} \times \mathbb{R}^\infty$ to $\mathcal{B}$, $F_{b_n}$ being the function on $[0, \tau] \in \mathcal{B}$ associated with $(n, b_n)$ defined in (1.1), $X_1, \ldots, X_K$ being constant design points, $\mathcal{F}$ being the Borel $\sigma$-field on $\mathcal{B} \times \mathbb{R}^\infty$. We also assume the marginal distribution of $\mathcal{P}$ on $\mathcal{B}$ is the prior $\pi$.

Our purpose is to illustrate a Bayesian regression method and the above mathematical formulation is only meant to facilitate an simple formal presentation. In fact, $\mathcal{P}$ is constructed after the prior and the likelihood are specified. A natural way to introduce the prior $\pi$ is to define $\pi(n, b_n) = p(n)\pi_n(b_n)$, with $\sum_{n=1}^{\infty} p(n) = 1$ and $\pi_n$ a density function on $B_n$; in fact, $\pi_n(\cdot)$ is the conditional density of $\pi$ on $B_n$ and also denoted by $\pi(\cdot \mid \{n\} \times B_n)$. Given $\mathbf{B} = (n, b_n)$, the likelihood for the data $\{(X_k, Y_{jk}) \mid j = 1, \ldots, m_k; k = 1, \ldots, K\}$ is

$$\prod_{k=1}^{K} \prod_{j=1}^{m_k} g_k(Y_{jk} - F_{b_n}(X_k)).$$

Thus the posterior density $\nu$ of the parameter $(n, b_n)$ given the data is proportional to

$$\prod_{k=1}^{K} \prod_{j=1}^{m_k} g_k(Y_{jk} - F_{b_n}(X_k))\pi_n(b_n)p(n),$$

where $(n, b_n) \in \mathcal{B}$. We note that, although we assume $g_k$ is known in this paper, the method of this paper can be extended to treat the case that $g_k$ has certain parametric form with priors on the parameters.

### 2.3. Support of Bernstein priors

The following propositions show that the support of the Bernstein priors can be quite large.

**Proposition 3.** *(Monotone) Let $B_n^{(1)} = \{b_n \in B_n : F_{b_n} \in I_n^{(1)}\}$. Assume $p(n) > 0$ for $n = 1, 2, \ldots$, and the conditional density $\pi_n(b_{0,n}, b_{1,n}, \ldots, b_{n,n})$ of $\pi(\cdot \mid \{n\} \times B_n^{(1)})$ has support $B_n^{(1)}$ for infinitely often $n$. Let $F$ be a given increasing and continuous function on $[0, \tau]$. Then $\pi\{(n, b_n) \in \bigcup_{n=2}^{\infty}(\{n\} \times B_n^{(1)}) : \|F_{b_n} - F\|_\infty < \epsilon\} > 0$ for every $\epsilon > 0$.*

*Proof.* Let $F_{(n)}(t) = \sum_{i=0}^{n} F(i\tau/n)\varphi_{i,n}(t/\tau)$. Using Bernstein-Weierstrass Theorem, we can choose a large $n_1$ so that $\|F_{(n_1)} - F\|_\infty \leq \epsilon/2, P(n_1) > 0$, and $\pi(\cdot \mid \{n_1\} \times B_{n_1}^{(1)}) > 0$ has support $B_{n_1}^{(1)}$. Combining this with $\|F_{b_n} - F_{(n)}\|_\infty \leq \max_{i=0,\ldots,n}\{|b_{i,n} - F(i\tau/n)|\}$ for $b_n \in B_n$, we get

$$\pi\{(n, b_n) \in \mathcal{B} : \|F_{b_n} - F\|_\infty < \epsilon\}$$
$$\geq \pi\{(n_1, b_{n_1}) \in \mathcal{B} : \|F_{b_{n_1}} - F_{(n_1)}\|_\infty < \frac{\epsilon}{2}\}$$
$$\geq \pi\{(n_1, b_{n_1}) \in \{n_1\} \times B_{n_1}^{(1)} : \max_{i=0,\ldots,n_1}|b_{i,n_1} - F(\frac{i\tau}{n_1})| < \frac{\epsilon}{2}\},$$

which is positive. This completes the proof. □

**Remarks.** If we know $c_1 < F(\tau) < c_0$, then it suffices to assume $\pi_n$ has support $\{(b_{0,n}, \ldots, b_{n,n}) \in B_n^{(1)} \mid c_1 \leq b_{0,n}, b_{n,n} \leq c_0\}$ in Proposition 3. Statements similar to Proposition 3 can also be made for concave functions. In fact, we have

**Proposition 4.** *(Unimodal Concave) Let $B_n^{(2)} = \{b_n \in B_n : F_{b_n} \in I_n^{(2)}\}$. Assume $p(n) > 0$ for $n = 1, 2, \ldots$, and the conditional density $\pi_n(b_{0,n}, b_{1,n}, \ldots, b_{n,n})$ of $\pi(\cdot \mid \{n\} \times B_n^{(2)})$ has support $B_n^{(2)}$ for infinitely often $n$. Let $F$ be a continuously differentiable real-valued function defined on $[0, \tau]$ with $F'(0) \geq 0$, $F'(\tau) \leq 0$, and $F'$ decreasing. Then $\pi\{(n, b_n) \in \bigcup_{n=2}^{\infty}(\{n\} \times B_n^{(2)}) : \|F_{b_n} - F\|_\infty + \|F'_{b_n} - F'\|_\infty < \epsilon\} > 0$ for every $\epsilon > 0$.*

**Proposition 5.** *(Unimodal) Let $B_n^{(3)} = \{b_n \in B_n : F_{b_n} \in I_n^{(3)}\}$. Assume $p(n) > 0$ for $n = 1, 2, \ldots$, and the conditional density $\pi_n(b_{0,n}, b_{1,n}, \ldots, b_{n,n})$ of $\pi(\cdot \mid \{n\} \times B_n^{(3)})$ has support $B_n^{(3)}$ for infinitely often $n$. Let $F$ be a continuously differentiable real-valued function defined on $[0, \tau]$ satisfying the properties that $F'(0) = 0$ and that there exists $s \in (0, \tau)$ such that $F'(s) = 0, F'(x) \geq 0$ for $x \in [0, s]$, and $F'(x) \leq 0$ for $x \in [s, \tau]$. Then $\pi\{(n, b_n) \in \bigcup_{n=2}^{\infty}(\{n\} \times B_n^{(3)}) : \|F_{b_n} - F\|_\infty + \|F'_{b_n} - F'\|_\infty < \epsilon\} > 0$ for every $\epsilon > 0$.*

## 3. Bayesian inference

Instead of defining priors explicitly, we propose algorithms to generate samples from the Bernstein priors that incorporate geometric information. We also propose algorithms for generating posterior distributions so as to do statistical inference. We consider only isotonic regression and unimodal concave (convex) regression in the rest of this paper, because we want to compare our method with the density-regression method. Algorithms 1, 2 and 3 respectively generate prior for isotonic regression, unimodal concave regression and unimodal convex regression. Both Algorithm 4, an independent Metropolis algorithm (IMA), and Algorithm 5, a Metropolis-Hastings reversible jump algorithm (MHRA), can be used to generate posterior for isotonic regression. Algorithms generating posterior for unimodal concave (convex) regression can be similarly proposed; they are omitted to make the paper more concise.

**Algorithm 1.** (*Bernstein isotonic prior*)

Let $p(1) = e^{-\alpha} + \alpha e^{-\alpha}$, $p(n) = \alpha^n e^{-\alpha}/n!$ for $n = 2, \ldots, n_0 - 1$, and $p(n_0) = 1 - \sum_{n=1}^{n_0-1} p(n)$. Let $q_1$ be a density so that its support contains $F(0)$. Let $q_2$ be a density so that its support contains $F(\tau)$. The following steps provide a way to sample from an implicitly defined prior distribution.

1. Generate $n \sim p$.
2. Generate $a_0 \sim q_1$, $a_n \sim q_2$.
3. Let $U_1, U_2, \ldots, U_{n-1}$ be a random sample from $Uniform(a_0, a_n)$. Let $U_{(1)}$, $U_{(2)}, \ldots, U_{(n-1)}$ be the order statistics of $\{U_1, U_2, \ldots, U_{n-1}\}$. Set $a_1 = U_{(1)}$, $a_2 = U_{(2)}, \ldots, a_{n-1} = U_{(n-1)}$.
4. $(n, a_0, \ldots, a_n) \in \mathcal{B}$ is a sample from the prior.

**Algorithm 2.** (*Bernstein concave prior*)

Let $p(2) = (1 + \alpha)e^{-\alpha} + \alpha^2 e^{-\alpha}/2$, $p(n) = \alpha^n e^{-\alpha}/n!$ for $n = 3, \ldots, n_0 - 1$, and $p(n_0) = 1 - \sum_{n=2}^{n_0-1} p(n)$. Let $q$ be a density with its support containing the mode of the regression function. Let $\beta_1$ be a lower bound of $F(0)$. Let $\beta_2$ be a lower bound of $F(\tau)$. The algorithm has the following steps.

1. Generate $n \sim p$.
2. Randomly choose $l$ from $\{1, 2, \ldots, n - 1\}$.
3. Generate $a_l \sim q$.
4. Generate $a_0 \sim Uniform(-a_l + 2\beta_1, a_l)$ and $a_n \sim Uniform(-a_l + 2\beta_2, a_l)$.
5. Let $U_1 \leq U_2 \leq \cdots \leq U_{l-1}$ be the order statistics of a random sample from $Uniform(a_0, a_l)$. Denote by $c_1 \leq c_2 \leq \cdots \leq c_l$ the order statistics of $\{a_l - U_{l-1}, U_{l-1} - U_{l-2}, \ldots, U_2 - U_1, U_1 - a_0\}$. Then set $a_1 = a_0 + c_l, a_2 = a_0 + c_l + c_{l-1}, \ldots, a_{l-1} = a_0 + c_l + \cdots + c_2$.
6. Let $V_1 \leq V_2 \leq \cdots \leq V_{n-l-1}$ be the order statistics of a random sample from $Uniform(a_n, a_l)$. Denote by $c'_1 \leq c'_2 \leq \cdots \leq c'_{n-l}$ the order statistics of $\{a_l - V_{n-l-1}, V_{n-l-1} - V_{n-l-2}, \ldots, V_2 - V_1, V_1 - a_n\}$. Then set $a_{l+1} = a_l - c'_1, a_{l+2} = a_l - c'_1 - c'_2, \ldots, a_{n-1} = a_l - c'_1 - \cdots - c'_{n-l-1}$.

In the above algorithms, the conditional distributions of $\pi(\cdot \mid \{n\} \times B_n)$ are defined to be those of $(a_0, a_1, \ldots, a_n)$. Although these two algorithms might look a little ad hoc, the main idea is to produce random sequence $a_0, a_1, \ldots, a_n$ satisfying conditions in the propositions in Section 2. It is obvious that $a_0 \leq a_1 \leq \cdots \leq a_n$ in Algorithm 1 and that $a_1 - a_0 > 0$, $a_n - a_{n-1} < 0$, and $a_{i+2} + a_i \leq 2a_{i+1}$ in Algorithm 2. It follows from Proposition 1 that the support of the prior generated by Algorithm 1 (Algorithm 2) contains only monotone (unimodal concave) functions. The following Algorithm 3 will be used in the simulation study.

**Algorithm 3.** (*Bernstein convex prior*)

Let $p(2) = (1 + \alpha)e^{-\alpha} + \alpha^2 e^{-\alpha}/2$, $p(n) = \alpha^n e^{-\alpha}/n!$ for $n = 3, \ldots, n_0 - 1$, and $p(n_0) = 1 - \sum_{n=2}^{n_0-1} p(n)$. Let $q$ be a density with its support containing the minimum value of the regression function. Let $\beta_1$ be a upper bound of $F(0)$, and $\beta_2$ be a upper bound of $F(\tau)$. The algorithm has the following steps.

1. Generate $n \sim p$.
2. Randomly choose $l$ from $\{1, 2, \ldots, n-1\}$.
3. Generate $a_l \sim q$.
4. Generate $a_0 \sim Uniform(a_l, 2\beta_1 - a_l)$ and $a_n \sim Uniform(a_l, 2\beta_2 - a_l)$.
5. Let $U_1 \leq U_2 \leq \cdots \leq U_{l-1}$ be the order statistics of a random sample from $Uniform(a_l, a_0)$. Denote by $c_1 \leq c_2 \leq \cdots \leq c_l$ the order statistics of $\{a_0 - U_{l-1}, U_{l-1} - U_{l-2}, \ldots, U_2 - U_1, U_1 - a_l\}$. Then set $a_1 = a_0 - c_l, a_2 = a_0 - c_l - c_{l-1}, \ldots, a_{l-1} = a_0 - c_l - \cdots - c_2$.
6. Let $V_1 \leq V_2 \leq \cdots \leq V_{n-l-1}$ be the order statistics of a random sample from $Uniform(a_l, a_n)$. Denote by $c_1' \leq c_2' \leq \cdots \leq c_{n-l}'$ the order statistics of $\{a_n - V_{n-l-1}, V_{n-l-1} - V_{n-l-2}, \ldots, V_2 - V_1, V_1 - a_l\}$. Then set $a_{l+1} = a_l + c_1', a_{l+2} = a_l + c_1' + c_2', \ldots, a_{n-1} = a_l + c_1' + \cdots + c_{n-l-1}'$.

**Algorithm 4.** (*IMA for the posterior in isotonic regression*)

This algorithm uses the independent Metropolis approach (see, for example, Robert and Casella [29], page 276) to calculate the posterior distribution $\nu$ of $(n, b)$. Let $x = (n, a_0, \ldots, a_n)$ be generated by the prior. We describe the transition from the current state $x^{(t)} = (n', a_0', \ldots, a_{n'}')$ to a new point $x^{(t+1)}$ by

$$x^{(t+1)} = \begin{cases} x, & \text{with prob. } \min\left\{1, \dfrac{\nu(x)\pi_{n'}(a_0', \ldots, a_{n'}')p(n')}{\nu(x^{(t)})\pi_n(a_0, \ldots, a_n)p(n)}\right\}, \\ x^{(t)}, & \text{o.w.} \end{cases}$$

The posterior distribution $\nu$ of $(n, a)$ in turn produces the posterior distribution of $F_a$, and the Bayes estimate we consider is the mean of $F_a$.

**Algorithm 5.** (*MHRA for the posterior in isotonic regression*)

Let $B_{(n)} = \{(n, a_0, \ldots, a_n) \,|\, (a_0, \ldots, a_n) \in B_n^{(1)}\}$. Let $x^{(t)} = (n, a_0, \ldots, a_n) \in B_{(n)}$ be the current state. We describe the transition from $x^{(t)} \in B_{(n)}$ to a new point $x^{(t+1)}$ by the following algorithms.

Randomly select one of three types of moves $H$, $H^+$, or $H^-$. Here $H$ is a transition of element in $B_{(n)}$, $H^+$ a transition of element from $B_{(n)}$ to $B_{(n+1)}$, and $H^-$ a transition of element from $B_{(n)}$ to $B_{(n-1)}$, respectively. Denote by $P_H^n$, $P_{H^+}^n$ and $P_{H^-}^n$ the probabilities of selecting the three different types of moves $H$, $H^+$ and $H^-$ when the current state of the Markov chain is in $B_{(n)}$. We set $P_{H^-}^1 = P_{H^+}^{n_0} = 0$. Consider $P_H^n = 1 - P_{H^+}^n - P_{H^-}^n$, $P_{H^+}^n = c \min\{1, \frac{p(n+1)}{p(n)}\}$ and $P_{H^-}^n = c \min\{1, \frac{p(n-1)}{p(n)}\}$, where $p$ is given in Algorithm 1 and $c$ is a sample parameter. Suppose $M_1 \leq F \leq M_2$.

If the move of type $H$ is selected, then

1. select $k$ randomly from $\{0, 1, 2, \ldots, n\}$ so that there is $1/3$ probability of choosing 0 or $n$ and $1/3(n-1)$ probability of choosing any one of the rest; generate $V \sim Uniform(c_{k-1}, c_{k+1})$, with $c_{-1} = M_1, c_{n+1} = M_2$, and $c_k = a_k$ for $k = 1, \ldots, n$;
2. let $y^{(t)}$ be the vector $x^{(t)}$ with $a_k$ replaced by $V$;
3. set the next state

$$x^{(t+1)} = \begin{cases} y^{(t)}, & \text{with prob. } \rho = \min\{1, \dfrac{\nu(y^{(t)})}{\nu(x^{(t)})}\}, \\ x^{(t)}, & \text{o.w.} \end{cases}$$

If the move of type $H^+$ is selected, then

1. generate $V \sim Uniform(a_0, a_n)$ and assume $a_k \leq V < a_{k+1}$;
2. let $y^{(t)} = (n+1, a_0, a_1, \ldots, a_k, V, a_{k+1}, \ldots, a_n) \in B_{(n+1)}$;
3. set the next state

$$x^{(t+1)} = \begin{cases} y^{(t)}, & \text{with prob. } \rho = \min\{1, \dfrac{\nu(y^{(t)}) \times p(n) \times (a_n - a_0)}{\nu(x^{(t)}) \times p(n+1) \times (n)}\}, \\ x^{(t)}, & \text{o.w.} \end{cases}$$

If the move of type $H^-$ is selected, then

1. select $k$ uniformly from $\{1, 2, \ldots, n-1\}$;
2. let $y^{(t)} = (n-1, a_0, a_1, \ldots, a_{k-1}, a_{k+1}, \ldots, a_n) \in B_{(n-1)}$;
3. set the next state

$$x^{(t+1)} = \begin{cases} y^{(t)}, & \text{with prob. } \rho = \min\{1, \dfrac{\nu(y^{(t)}) \times p(n) \times (n-1)}{\nu(x^{(t)}) \times p(n-1) \times (a_n - a_0)}\}, \\ x^{(t)}, & \text{o.w.} \end{cases}$$

## 4. Numerical studies

We now explore the numerical performance of the Bayes methods in this paper. Viewing the posterior as a distribution on regression functions, we can use the posterior mean $\hat{F}$ of the regression functions as the estimate; namely,

$$\hat{F}(t) = \frac{1}{m} \sum_{j=1}^{m} F_{b^{(j)}}(t),$$

where $m$ is a large number, and $b^{(1)}, b^{(2)}, \ldots$ are chosen randomly according to the posterior distribution. The performance of $\hat{F}$ is evaluated by the $L_1$-norm, sup-norm and mean square error (MSE) of $\hat{F} - F$, as a function on $[0, 1]$. Here $F$ denotes the true regression function. We note that, in the studies of isotonic regression and concave (convex) regression, $\hat{F}$ is a reasonable estimate because monotonicity and concavity (convexity) are preserved by linear combination with non-negative coefficients, and when studying other shape restricted regression, it might be more appropriate to use posterior mode as the estimate.

Numerical studies in this section include comparison between our method and the density-regression method. In order to make the comparison more convincing, we study in this section both isotonic regression and convex regression, instead of concave regression, because they are studied by Dette and coauthors.

In this section, $\tau = 1$; $K = 100$; $X_1, X_2, \ldots, X_{100}$ are i.i.d. from $Uniform(0, 1)$; $m_k = 1$ for every $k = 1, 2, \ldots, 100$; $g_k$ is $Normal(0, \sigma^2)$ with $\sigma = 0.1$ or $1$ in the data generation. When carrying out inference, $\sigma$ is estimated from the data.

### 4.1. Isotonic regression

Our simulation studies suggest that compared with the density-regression method of Dette et al. [11], our method performs comparably when the noise is large and better when the noise is small. We studied all the regression functions in Dette et al. [11], Dette and Pilz [12], and Dette [10], and all the results are similar; hence, we only report the results for two of the regression functions, which are defined by

$$F_{<1>}(t) = \sin(\frac{\pi}{2}t), \tag{4.1}$$

$$F_{<2>}(t) = \begin{cases} 2t & \text{for } t \in [0, 0.25], \\ 0.5 & \text{for } t \in [0.25, 0.75], \\ 2t - 1 & \text{for } t \in [0.75, 1]. \end{cases} \qquad (4.2)$$

We note that these are not polynomials and our method usually performs better for polynomials. We also note that the notation $\Pi$ in (4.1) is the ratio of the circumference of a circle to its diameter.

We use Algorithm 1 with $n_0 = 20$, $\alpha = 10$, $q_1 \sim Uniform(q_{11}, q_{12})$, $q_2 \sim Uniform(q_{21}, q_{22})$ to generate the prior distribution. Here $q_{11}$, $q_{12}$, $q_{21}$ and $q_{22}$ are defined as follows. Let $X_{(1)} < X_{(2)} < \cdots < X_{(100)}$ be the order statistics for $\{X_1, X_2, \ldots, X_{100}\}$. Let $Y_{[i]} = Y_j$ if $X_{(i)} = X_j$. Then $q_{11} = \min\{Y_{[i]} \mid i = 1, 2, \ldots, 10\}$, $q_{12} = q_{21} = \frac{1}{100} \sum_{i=1}^{100} Y_i$ and $q_{22} = \max\{Y_{[i]} \mid i = 91, 92, \ldots, 100\}$. We note that the choice of $n_0 = 20$ and $\alpha = 10$ makes it relatively uninformative in the order of the polynomial.

We use Algorithm 5 (MHRA) with $c = 0.35$, $M_1 = q_{11}$, $M_2 = q_{22}$ and estimated $\sigma^2$ to generate the posterior distribution. We note that this choice of $c$ allows relatively large probabilities of changing the order of the polynomial and $\sigma^2$ is estimated by $\frac{1}{198} \sum_{i=1}^{99} (Y_{[i+1]} - Y_{[i]})^2$, which is also used in the following density-regression method. We run MHRA for 100,000 updates; after a burn in period of 10,000 updates, we use the remaining 90,000 realizations of the Markov chain to obtain the posterior mean. Starting from generating the data $\{X_i, Y_i \mid i = 1, 2, \ldots, 100\}$, the above experiment is replicated 200 times; the averages of the resulting 200 $L_1$-norms, sup-norms and MSE of $\hat{F} - F$ are reported in Table 1 and Table 2, which also contain the corresponding results using the density-regression method; the method of this paper is referred to as the Bayes method. Also contained in Table 1 and Table 2 are the posterior distributions of the polynomial order; these posterior distributions seem to suggest that $n_0 = 20$ in the prior is large enough. Figure 1 contains the autocorrelation plots of the $L_1$-norm of $F_{b^{(j)}}$ for one of the 200 replicates in the study of $F_{<1>}$. Figure 1 indicates that MHRA behaves quite nicely. Those for $F_{<2>}$ look similar and hence are omitted. The corresponding results using IMA are also omitted, because they are similar to those using MHRA, except having larger autocorrelation values and smaller effective sample sizes. Concepts of autocorrelation and effective sample size can be found in Liu [19], for example.

### 4.2. Convex regression

The main conclusion of our simulation studies regarding convex regression is similar to that reported for isotonic regression; namely, compared with the density-regression method for convex regression studied by Birke and Dette [3], our Bayes method performs comparably when the noise is large and performs better when the noise is small. To make the presentation concise, we only report the results for two of the regression functions, which are defined by

$$F_{<3>}(t) = (16/9)(t - 1/4)^2, \qquad (4.3)$$

$$F_{<4>}(t) = \begin{cases} -4t + 1 & \text{for } t \in [0, 0.25], \\ 0 & \text{for } t \in [0.25, 0.75], \\ 4t - 3 & \text{for } t \in [0.75, 1]. \end{cases} \qquad (4.4)$$

We note that both $F_{<3>}$ and $F_{<4>}$ are examples in Birke and Dette [3].

TABLE 1
*Simulation study for $F_{<1>}(t) = \sin(\frac{\pi}{2}t)$.*

| | | a | |
|---|---|---|---|
| $\sigma$ | $|\hat{F} - F|$ | Bayes | Density -regression |
| | $L_1$-norm | 0.0161 | 0.0238 |
| 0.1 | Sup-norm | 0.0532 | 0.0874 |
| | MSE | 0.0004 | 0.0010 |
| | $L_1$-norm | 0.1148 | 0.1194 |
| 1 | Sup-norm | 0.2795 | 0.2651 |
| | MSE | 0.0215 | 0.0233 |

| | | | | | | | | b |
|---|---|---|---|---|---|---|---|---|
| $\sigma$ | $n$ | Posterior probability | $n$ | Posterior probability | $n$ | Posterior probability | $n$ | Posterior probability |
| | 1 | 0.0000 | 6 | 0.0797 | 11 | 0.0960 | 16 | 0.0139 |
| | 2 | 0.0300 | 7 | 0.0994 | 12 | 0.0708 | 17 | 0.0076 |
| 0.1 | 3 | 0.0433 | 8 | 0.1009 | 13 | 0.0497 | 18 | 0.0037 |
| | 4 | 0.0483 | 9 | 0.1109 | 14 | 0.0359 | 19 | 0.0016 |
| | 5 | 0.0789 | 10 | 0.1044 | 15 | 0.0235 | 20 | 0.0015 |
| | 1 | 0.0006 | 6 | 0.0673 | 11 | 0.1104 | 16 | 0.0200 |
| | 2 | 0.0026 | 7 | 0.0950 | 12 | 0.0878 | 17 | 0.0113 |
| 1 | 3 | 0.0077 | 8 | 0.1164 | 13 | 0.0687 | 18 | 0.0068 |
| | 4 | 0.0196 | 9 | 0.1308 | 14 | 0.0495 | 19 | 0.0042 |
| | 5 | 0.0406 | 10 | 0.1252 | 15 | 0.0319 | 20 | 0.0036 |

We use Algorithm 3 with $n_0 = 20$, $\alpha = 10$, $q \sim Uniform(q_{01}, q_{02})$, $\beta_1 = q_3$ and $\beta_2 = q_4$ to generate the prior distribution, where $q_{01}$, $q_{02}$, $q_3$ and $q_4$ are defined as follows. Let $Y_{(1)} < Y_{(2)} < \cdots < Y_{(100)}$ be the order statistics for $\{Y_1, Y_2, \ldots, Y_{100}\}$. Then $q_{01} = \frac{1}{10}\sum_{i=1}^{10} Y_{(i)}$, $q_{02} = |q_{01} + \frac{1}{100}\sum_{i=1}^{100} Y_i|/2$. Let $X_{(1)} < X_{(2)} < \cdots < X_{(100)}$ be the order statistics for $\{X_1, X_2, \ldots, X_{100}\}$. Let $Y_{[i]} = Y_j$ if $X_{(i)} = X_j$. Then $q_3 = \max\{Y_{[i]} \mid i = 1, 2, \ldots, 5\}$ and $q_4 = \max\{Y_{[i]} \mid i = 96, 97, \ldots, 100\}$.

We use MHRA to generate the posterior distribution with similar parameters given in Algorithm 5. Table 3 and Table 4 contain the main results. Entries in Table 3 and Table 4 bear the same meanings of the corresponding entries in Table 1. Figure 2 carries similar information as Figure 1.

Table 1a gives comparison between the Bayes method and the density-regression method. Table 1b gives the posterior probability of the order of the Bernstein polynomial. The effective sample size for $\sigma = 0.1$ (1) is 5623 (461). The acceptance rate for $\sigma = 0.1$ (1) is 0.2923 (0.6635).

Table 2a gives comparison between the Bayes method and the density-regression method. Table 2b gives the posterior probability of the order of the Bernstein polynomial. The effective sample size for $\sigma = 0.1$ (1) is 4913 (364). The acceptance rate for $\sigma = 0.1$ (1) is 0.3779 (0.6943).

Table 3a gives comparison between the Bayes method and the density-regression method. Table 3b gives the posterior probability of the order of the Bernstein polynomial. The effective sample size for $\sigma = 0.1$ (1) is 1397 (344). The acceptance rate for $\sigma = 0.1$ (1) is 0.2413 (0.3430).

Table 4a gives comparison between the Bayes method and the density-regression method. Table 4b gives the posterior probability of the order of the Bernstein polynomial. The effective sample size for $\sigma = 0.1$ (1) is 722 (216). The acceptance rate for $\sigma = 0.1$ (1) is 0.2217 (0.3670).

TABLE 2

*Simulation study for* $F_{<2>}(t) = \begin{cases} 2t & \text{for } t \in [0, 0.25], \\ 0.5 & \text{for } t \in [0.25, 0.75], \\ 2t - 1 & \text{for } t \in [0.75, 1]. \end{cases}$

a

| $\sigma$ | $|\hat{F} - F|$ | Bayes | Density -regression |
|---|---|---|---|
| | $L_1$-norm | 0.0353 | 0.0413 |
| 0.1 | Sup-norm | 0.1002 | 0.1415 |
| | MSE | 0.0019 | 0.0030 |
| | $L_1$-norm | 0.1255 | 0.1267 |
| 1 | Sup-norm | 0.3035 | 0.3509 |
| | MSE | 0.0244 | 0.0256 |

b

| $\sigma$ | $n$ | Posterior probability | $n$ | Posterior probability | $n$ | Posterior probability | $n$ | Posterior probability |
|---|---|---|---|---|---|---|---|---|
| | 1 | 0.0000 | 6 | 0.0144 | 11 | 0.1481 | 16 | 0.0375 |
| | 2 | 0.0000 | 7 | 0.0385 | 12 | 0.1730 | 17 | 0.0194 |
| 0.1 | 3 | 0.0000 | 8 | 0.0420 | 13 | 0.1571 | 18 | 0.0117 |
| | 4 | 0.0000 | 9 | 0.0621 | 14 | 0.1142 | 19 | 0.0062 |
| | 5 | 0.0005 | 10 | 0.1006 | 15 | 0.0692 | 20 | 0.0055 |
| | 1 | 0.0004 | 6 | 0.0611 | 11 | 0.1169 | 16 | 0.0195 |
| | 2 | 0.0016 | 7 | 0.0871 | 12 | 0.0992 | 17 | 0.0113 |
| 1 | 3 | 0.0062 | 8 | 0.1100 | 13 | 0.0794 | 18 | 0.0060 |
| | 4 | 0.0177 | 9 | 0.1245 | 14 | 0.0539 | 19 | 0.0029 |
| | 5 | 0.0379 | 10 | 0.1284 | 15 | 0.0335 | 20 | 0.0025 |

TABLE 3
*Simulation study for* $F_{<3>}(t) = (16/9)(t - 1/4)^2$.

a

| $\sigma$ | $|\hat{F} - F|$ | Bayes | Density -regression |
|---|---|---|---|
| | $L_1$-norm | 0.0157 | 0.0775 |
| 0.1 | Sup-norm | 0.0525 | 0.2296 |
| | MSE | 0.0004 | 0.0101 |
| | $L_1$-norm | 0.1362 | 0.1366 |
| 1 | Sup-norm | 0.4247 | 0.3643 |
| | MSE | 0.0319 | 0.0292 |

b

| $\sigma$ | $n$ | Posterior probability | $n$ | Posterior probability | $n$ | Posterior probability | $n$ | Posterior probability |
|---|---|---|---|---|---|---|---|---|
| | 2 | 0.0001 | 7 | 0.1177 | 12 | 0.0452 | 17 | 0.0064 |
| | 3 | 0.0609 | 8 | 0.1312 | 13 | 0.0346 | 18 | 0.0029 |
| 0.1 | 4 | 0.0826 | 9 | 0.1139 | 14 | 0.0230 | 19 | 0.0020 |
| | 5 | 0.0931 | 10 | 0.0938 | 15 | 0.0147 | 20 | 0.0021 |
| | 6 | 0.1015 | 11 | 0.0653 | 16 | 0.0090 | | |
| | 2 | 0.0413 | 7 | 0.1094 | 12 | 0.0516 | 17 | 0.0062 |
| | 3 | 0.0602 | 8 | 0.1175 | 13 | 0.0313 | 18 | 0.0028 |
| 1 | 4 | 0.0728 | 9 | 0.1110 | 14 | 0.0204 | 19 | 0.0015 |
| | 5 | 0.0849 | 10 | 0.0886 | 15 | 0.0141 | 20 | 0.0014 |
| | 6 | 0.1037 | 11 | 0.0717 | 16 | 0.0096 | | |

TABLE 4

*Simulation study for* $F_{<4>}(t) = \begin{cases} -4t+1 & \text{for } t \in [0, 0.25], \\ 0 & \text{for } t \in [0.25, 0.75], \\ 4t-3 & \text{for } t \in [0.75, 1]. \end{cases}$

a

| $\sigma$ | $|\hat{F} - F|$ | Bayes | Density -regression |
|---|---|---|---|
| 0.1 | $L_1$-norm | 0.0405 | 0.1311 |
| | Sup-norm | 0.1381 | 0.4663 |
| | MSE | 0.0026 | 0.0282 |
| 1 | $L_1$-norm | 0.1394 | 0.2265 |
| | Sup-norm | 0.4603 | 0.7010 |
| | MSE | 0.0338 | 0.0793 |

b

| $\sigma$ | $n$ | Posterior probability | $n$ | Posterior probability | $n$ | Posterior probability | $n$ | Posterior probability |
|---|---|---|---|---|---|---|---|---|
| | 2 | 0.0000 | 7 | 0.0000 | 12 | 0.1710 | 17 | 0.0263 |
| | 3 | 0.0000 | 8 | 0.0205 | 13 | 0.1301 | 18 | 0.0158 |
| 0.1 | 4 | 0.0000 | 9 | 0.1144 | 14 | 0.1130 | 19 | 0.0124 |
| | 5 | 0.0000 | 10 | 0.1028 | 15 | 0.0778 | 20 | 0.0187 |
| | 6 | 0.0000 | 11 | 0.1498 | 16 | 0.0474 | | |
| | 2 | 0.0153 | 7 | 0.1194 | 12 | 0.0571 | 17 | 0.0068 |
| | 3 | 0.0295 | 8 | 0.1254 | 13 | 0.0432 | 18 | 0.0041 |
| 1 | 4 | 0.0539 | 9 | 0.1153 | 14 | 0.0321 | 19 | 0.0013 |
| | 5 | 0.0784 | 10 | 0.0953 | 15 | 0.0227 | 20 | 0.0014 |
| | 6 | 0.1065 | 11 | 0.0786 | 16 | 0.0137 | | |

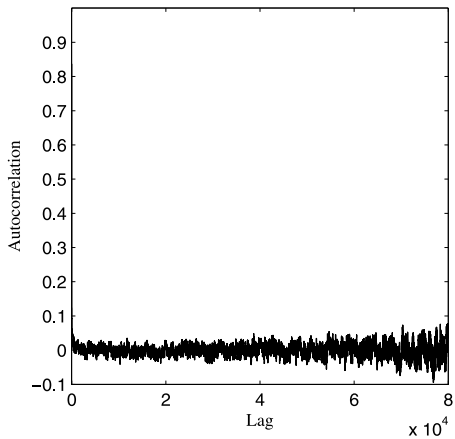
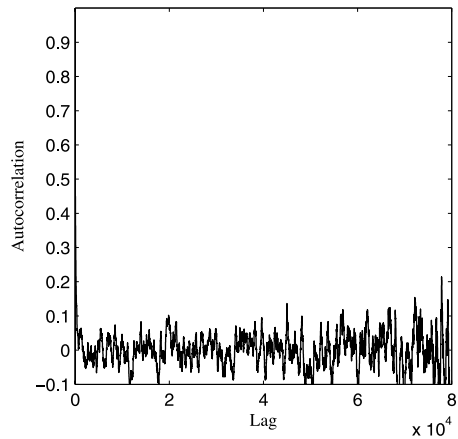

Figure 1a: σ = 0.1          Figure 1b: σ = 1

FIG 1. *Autocorrelation plots for the MHRA in the data generation from the posterior distribution for* $F_{<1>}(t) = \sin(\frac{\Pi}{2}t)$.

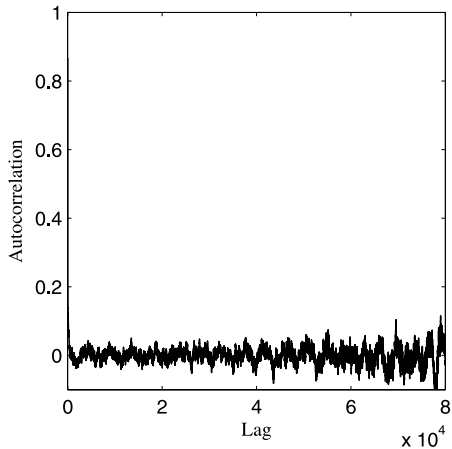

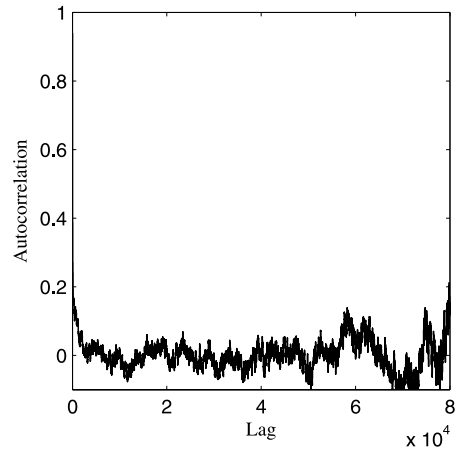

Figure 2a: $\sigma = 0.1$               Figure 2b: $\sigma = 1$

FIG 2. *Autocorrelation plots for the MHRA in the data generation from the posterior distribution for* $F_{<3>}(t) = (16/9)(t - 1/4)^2$.

## 5. Discussion

We have proposed a Bayes approach to shape restricted regression with the prior introduced by random Bernstein polynomials. In particular, algorithms for generating the priors and the posteriors are proposed and numerical performance of the Bayes estimate has been examined in some detail. The usefulness of this approach is successfully demonstrated in simulation studies for isotonic regression and convex regression, which compares our method with the density-regression method. We note that certain frequentist properties of this Bayes estimate are established in Chang et al. [6]. The method of this paper can also be used to assess the validity of the shape restriction on the regression function, by considering the predictive posterior assessment proposed by Gelman et al. [14] and Gelman [13]; in fact, we have implemented this assessment and found it quite satisfactory. We would like to remark that this approach can be used to study other shape restricted statistical inference problem. In fact, as pointed out by P. Bickel and M. Woodroofe in the Vardi memorial conference, the geometry of Bernstein polynomials may be utilized to propose a penalized likelihood approach to shape restricted regression; our preliminary simulation studies (not reported here) do suggest that penalized maximum likelihood estimate looks promising, and further investigation is underway.